\documentclass[12pt]{amsart}

\usepackage{amsthm, amssymb}
\usepackage{amsfonts, amscd}

\usepackage{epsfig,multicol}
\begin{document}
\title[Equivariant cohomology  and  analytic descriptions of ring isomorphisms]{\large \bf
 Equivariant cohomology and  analytic descriptions of ring isomorphisms }
\renewcommand{\thefootnote}{\fnsymbol{footnote}}
\author[Bo Chen and Zhi L\"u]{Bo Chen and Zhi L\"u}
\footnote[0]{ {\bf 2000 Mathematics Subject Classification.} 57S17,
55N91, 58J20
\endgraf
 {\bf Key words and phrases. } $G$-manifold, equivariant index, equivariant cohomology, analytic ring isomorphism.
\endgraf
  Supported by grants from NSFC (No. 10371020 and No. 10671034). }
\address{School of Mathematical Science, Fudan University, Shanghai,
200433, People's Republic of China.} \email{061018005@fudan.edu.cn}
\address{Institute of Mathematics, School of Mathematical Science, Fudan University, Shanghai,
200433, People's Republic of China.}
 \email{zlu@fudan.edu.cn}

\maketitle
\newtheorem{prop}{Proposition}[section]
\newtheorem{lem}{Lemma}[section]
\newtheorem{thm}{Theorem}[section]
\newtheorem{cor}{Corollary}[section]
\newtheorem{note}{Note}[section]

\theoremstyle{definition}

\newtheorem{defn}[thm]{Definition}
\newtheorem{notation}[thm]{Notation}
\newtheorem{example}[thm]{Example}
\newtheorem{conj}[thm]{Conjecture}
\newtheorem{prob}[thm]{Problem}

\theoremstyle{remark}

\newtheorem{rem}{Remark}

\def\mapleftright#1{\smash{\mathop{\longleftrightarrow}\limits^{#1}}}
\def \co{\colon\thinspace}
\def \R {\hbox{\rm I \kern -5.5pt R}}
\def \C {\Bbb C}
\def \Hom {{\rm Hom}}
\def \Im {{\rm Im}}
\def \dprm {{\prime\prime}}

\def\mapleftright#1{\smash{\mathop{\longleftrightarrow}\limits^{#1}}}
\def\mapleftright#1{\smash{\mathop{\longleftrightarrow}\limits^{#1}}}
\def \R {{\bf R}\/}
\def \Z {{\mathbb Z}\/}
\def \A {{\cal A}\/}
\def \G{{\rm GL}}
\def \vphi {{\varphi}\/}
\def \ds {\displaystyle}
\def\co{\colon\thinspace}
\def\a {\rm Aut}
\def\normalbaselines{\baselineskip20pt \lineskip3pt \lineskiplimit3pt }
\def\mapright#1{\smash{\mathop{\longrightarrow}\limits^{#1}}}
\def\mapdown#1{\Big\downarrow\rlap{$\vcenter{\hbox{$\scriptstyle#1$}}$}}
\def\mapup#1{\Big\uparrow\rlap{$\vcenter{\hbox{$\scriptstyle#1$}}$}}
\def\mapneup#1{\llap{$\vcenter{\hbox{$\scriptstyle#1$}}$}\nearrow}
\def\mapswdown#1{\swarrow\rlap{$\vcenter{\hbox{$\scriptstyle#1$}}$}}
\def\mapsedown#1{\llap{$\vcenter{\hbox{$\scriptstyle#1$}}$}\searrow}
\def\mapswdown#1{\swarrow\rlap{$\vcenter{\hbox{$\scriptstyle#1$}}$}}

\newcommand \cd {\mathbf{\cdot}}
\newcommand \mba {{\boldmath$\alpha$}}
\newcommand \mbe {\mathbf{e}\/}
\newcommand \mbf {\mathbf{f}\/}
\newcommand \mbg {\mathbf{g}\/}
\newcommand \mbh {\mathbf{h}\/}
\newcommand \mbi {\mathbf{i}\/}
\newcommand \mbj {\mathbf{j}\/}
\newcommand \mbk {\mathbf{k}\/}
\newcommand \bkp {{\bf\kappa}\/}
\newcommand \mbF {\mathbf{F}\/}
\newcommand \mbG {\mathbf{G}\/}
\newcommand \mbH {\mathbf{H}\/}
\newcommand \mbn {\mathbf{n}\/}
\newcommand \mbN {\mathbf{N}\/}
\newcommand \mbR {\mathbf{R}\/}
\newcommand \mbS {\mathbb{S}\/}
\newcommand \mbt {\mathbf{t}\/}
\newcommand \mbT {\mathbf{T}\/}
\newcommand \mbv {\mathbf{v}\/}
\newcommand \mbu {\mathbf{u}\/}
\newcommand \mvV {\mathbf{V}\/}
\newcommand \mbx {\mathbf{x}\/}
\newcommand \mby {\mathbf{y}\/}
\newcommand \mbz {\mathbf{z}\/}
\newcommand \mbvp {\mathbf{\varphi}\/}
\newcommand \dsp {\displaystyle}
\newcommand \vep {{\varepsilon}}
\newcommand \prm {\prime}

\addcontentsline{toc}{chapter}{Bibliography}

\begin{abstract} In this paper we consider a class of connected closed $G$-manifolds with a non-empty finite
fixed point set, each $M$ of which is  totally non-homologous to
zero in $M_G$ (or $G$-equivariantly formal), where $G={\Bbb Z}_2$.
With the help of the equivariant index, we give an explicit
description of the equivariant cohomology of such a $G$-manifold in
terms of algebra, so that we can  obtain analytic descriptions of
ring isomorphisms among equivariant cohomology rings of such
$G$-manifolds, and  a necessary and sufficient condition that the
equivariant cohomology rings of such two $G$-manifolds  are
isomorphic. This also leads us to analyze how many there are
equivariant cohomology rings up to isomorphism for such
$G$-manifolds in  2- and 3-dimensional cases.
\end{abstract}

\footnote[0]{}

\section{Introduction}

Throughout this paper,  assume that  $G={\Bbb Z}_2$ unless stated
otherwise. Let $EG\longrightarrow BG$ be the universal principal
$G$-bundle, where $BG=EG/G={\Bbb R}P^\infty$ is the classifying
space of $G$. It is well-known that $H^*(BG;{\Bbb Z}_2)={\Bbb
Z}_2[t]$ with the $t$ one-dimensional generator.

\vskip .2cm

Let $X$ be a $G$-space. Then $X_G:=EG\times_GX$-----the orbit space
of the diagonal action on the product $EG\times X$-----is the total
space of the bundle $X\longrightarrow X_G\longrightarrow BG$
associated to the universal principal bundle $G\longrightarrow
EG\longrightarrow BG$. The space $X_G=EG\times_GX$ is called the
{\em Borel construction} on the $G$-space $X$. Applying cohomology
with coefficients ${\Bbb Z}_2$ to  $X_G$  gives the {\em equivariant
cohomology} $H^*_G(X;{\Bbb Z}_2):=H^*(X_G;{\Bbb Z}_2).$ It is
well-known that equivariant cohomologies $H^*_G(X;{\Bbb Z}_2)$ and
$H^*_G(X^G;{\Bbb Z}_2)$ are $H^*(BG;{\Bbb Z}_2)$-modules; in
particular, $H^*_G(X^G;{\Bbb Z}_2)$ is a free $H^*(BG;{\Bbb
Z}_2)$-module.

\vskip .2cm

Suppose that $M$ is a connected closed manifold and admits a
$G$-action with $M^G$ a non-empty finite fixed set. For the
fibration $M\longrightarrow M_G\longrightarrow BG$, if the
restriction to a typical fiber $$H^*_G(M;{\Bbb Z}_2)\longrightarrow
H^*(M;{\Bbb Z}_2)$$ is an epimorphism, then $M$ is called {\em
totally non-homologous to zero} in $M_G$ (see \cite{b}). Under this
condition, $M$ is also called {\em $G$-equivariantly formal} (cf.
\cite{gkm}). Note that when the action group on $M$ is replaced by a
2-torus $({\Bbb Z}_2)^k$ with $k\geq 1$ and each component of $M^K$
has dimension at most 1 for $K< ({\Bbb Z}_2)^k$ a corank-1 2-torus,
if $M$ is $({\Bbb Z}_2)^k$-equivariantly formal (i.e., the edge
morphism $H^*_{({\Bbb Z}_2)^k}(M;{\Bbb Z}_2)\longrightarrow
H^*(M;{\Bbb Z}_2)$ is surjective), then there is a mod 2 GKM theory
(corresponding to GKM theory, see \cite{gkm}, \cite{gz}), indicating
that the $({\Bbb Z}_2)^k$-equivariant cohomology of $M$ can be
explicitly expressed in terms of its associated  graph $(\Gamma_M,
\alpha)$ (cf. \cite{bgh},  \cite{l}). In particular, when $k=1$, $M$
is naturally restricted to have dimension at most 1, so this means
that the GKM theory  can be carried out only for at most
 1-dimensional $G$-equivariantly formal manifolds.
 In this paper, we shall give  explicit descriptions of equivariant
cohomology rings of $G$-equivariantly formal manifolds at any
dimension, and these descriptions are of algebra rather than
combinatorics.

\vskip .2cm

Let $\Lambda_n$ denote the set of all $n$-dimensional connected
closed $G$-manifolds with a non-empty finite fixed point set, each
of which   is $G$-equivariantly formal. Note that obviously
$\Lambda_1$ contains a unique 1-dimensional closed manifold, i.e., a
circle $S^1$. Taking a $M$ in $\Lambda_n$, by Conner and Floyd
\cite{cf}, one knows that $|M^G|$ must be even. Let $r$ be a
positive integer, and write $\Lambda_n^{2r}=\big\{M\in\Lambda_n\big|
|M^G|=2r\big\}$. Then $\Lambda_n=\bigcup\limits_{r\geq
1}\Lambda_n^{2r}$. Given a $M$ in $\Lambda_n$, we know from
\cite{ap} and \cite{b} that the following conditions are equivalent

\begin{enumerate}
\item[(1)] $M$   is $G$-equivariantly formal;
 \item[(2)] $|M^G|=\sum\limits_{i=0}^nb_i$ where $b_i$ is the
 $i$-th mod 2 Betti number of $M$;
 \item[(3)] $H_G^*(M;{\Bbb Z}_2)$ is a free $H^*(BG;{\Bbb Z}_2)$-module;
\item[(4)] The inclusion $i:M^G\hookrightarrow M$
induces a monomorphism $i^*: H_G^*(M;{\Bbb Z}_2)\longrightarrow
H_G^*(M^G;{\Bbb Z}_2)$.
\end{enumerate}
If $|M^G|=2r$,  since $H^*_G(M^G;{\Bbb Z}_2)=\bigoplus\limits_{p\in
M^G}H^*_G(\{p\};{\Bbb Z}_2)$ and the equivariant cohomology of a
point is isomorphic to $H^*(BG;{\Bbb Z}_2)={\Bbb Z}_2[t]$, we have
that $H_G^*(M^G;{\Bbb Z}_2)\cong ({\Bbb Z}_2)^{2r}[t]$ is a
polynomial ring (or algebra). Thus we obtain a monomorphism from
$H_G^*(M;{\Bbb Z}_2)$ into $({\Bbb Z}_2)^{2r}[t]$, also denoted by
$i^*$, so $H_G^*(M;{\Bbb Z}_2)$ may be identified with a subring (or
subalgebra) of $({\Bbb Z}_2)^{2r}[t]$.

\vskip .2cm

Using the equivariant index, we shall give an explicit description
of $H_G^*(M; {\Bbb Z}_2)$ in $({\Bbb Z}_2)^{2r}[t]$ (see
Theorem~\ref{ring str}). Then we consider the following questions:

\begin{enumerate}
\item[(Q1)] If the equivariant cohomology rings of two $G$-manifolds
in $\Lambda_n$ are isomorphic, then can the isomorphism between them
be explicitly expressed?
 \item[(Q2)] When are the equivariant cohomology rings of two $G$-manifolds in $\Lambda_n
$ isomorphic?
 \item[(Q3)] How many are
there  equivariant cohomology rings (up to isomorphism) of
$G$-manifolds in $\Lambda_n$?
\end{enumerate}

We completely answer  (Q1) and (Q2). An interesting thing is that we
do not only find an explicit description for the isomorphism between
equivariant cohomology rings of two $G$-manifolds in $\Lambda_n$,
but such a description is also analytic (see Theorem~\ref{analyze}),
so that we may obtain a necessary and sufficient condition that the
equivariant cohomology rings of such two $G$-manifolds  are
isomorphic in terms of algebra (see Theorem~\ref{ns}).
 As for  (Q3), the question is answered completely in the case  $n=2$. When
$n=3$,  we find an upper bound of the number for the equivariant
cohomology rings (up to isomorphism) of  $G$-manifolds in
$\Lambda_3^{2r}$ (see Proposition~\ref{number}).

\vskip .2cm

The paper is organized as follows. In Section 2 we review the
localization theorem and reformulate the equivariant index from the
work of Allday and Puppe \cite{ap}. In Section 3 we study the
structure of equivariant cohomology of a $G$-manifold in $\Lambda_n$
and obtain an explicit description in terms of algebra. Then we
completely answer (Q2) and (Q3) in the  case $n=2$ in Section 4. In
Section 5, we give an analytic description for the isomorphism
between equivariant cohomology rings of two $G$-manifolds in
$\Lambda_n$, so that we may obtain a necessary and sufficient
condition that the equivariant cohomology rings of two $G$-manifolds
in $\Lambda_n$ are isomorphic. In section 6, we discuss the number
of the equivariant cohomology rings (up to isomorphism) of
$G$-manifolds in $\Lambda_3^{2r}$, and obtain an upper bound of the
number.

\vskip .2cm The authors  express their thanks to Shengzhi Xu for
helpful conversation in the argument of Lemma~\ref{bound}.
 The authors also
would like to express their gratitude to Professor Volker Puppe, who
informed of us that there is an essential relationship between the
equivariant cohomology rings and the  coding theory (see \cite{p}),
and Theorem 3 in \cite{kp}  implies that the map $g$ is surjective
in the remark~\ref{vp} of our paper.

\section{Localization theorem and equivariant index}

Suppose that $M$ is an $n$-dimensional  $G$-manifold with $M^G$ a
non-empty finite set. Let $S$ be the subset of $H^*(BG;{\Bbb Z}_2)$
 generated multiplicatively by nonzero elements in $H^1(BG;{\Bbb
 Z}_2)$. Then one has the following well-known localization theorem  (see \cite{ap}, \cite{h}).

 \begin{thm} [Localization theorem]\label{local}
  $$S^{-1}i^*:
 S^{-1}H^*_G(M;{\Bbb Z}_2)\longrightarrow S^{-1}H^*_G(M^G;{\Bbb
 Z}_2)$$
 is  an isomorphism of $S^{-1}H^*(BG;{\Bbb Z}_2)$-algebras, where $i$ is the inclusion of from
 $M^G$ into $M$.
 \end{thm}

Take an isolated point $p\in M^G$. Let $i_p$ be the inclusion of
from $p$ into $M$, then one has the equivariant Gysin homomorphism
$$i_{p!}: H^*_G(\{p\};{\Bbb Z}_2)\longrightarrow H^{*+n}_G(M;{\Bbb
Z}_2).$$ On the other hand, one has also  a natural induced
homomorphism
$$i^*_p: H^*_G(M;{\Bbb Z}_2)\longrightarrow H^*_G(\{p\};{\Bbb
Z}_2)$$ and in particular, it is easy to check that
$i^*=\bigoplus\limits_{p\in M^G}i_p^*$.  Furthermore, one knows that
the equivariant Euler class at $p$ is
$$\chi_G(p)=i^*_pi_{p!}(1_p)\in H^n_G(\{p\};{\Bbb Z}_2)= H^n(BG;{\Bbb Z}_2)={\Bbb Z}_2t^n,$$
which is equal to that of the real $G$-representation at $p$, where
$1_p\in H^*_G(\{p\};{\Bbb Z}_2)$ is the identity and ${\Bbb
Z}_2t^n=\{at^n|a\in{\Bbb Z}_2\}$. Thus, we may write
$\chi_G(p)=t^n$. Write $\theta_p=i_{p!}(1_p)$. Then $\theta_p\in
H^n_G(M;{\Bbb Z}_2)$ and $i^*_p(\theta_p)=\chi_G(p)$.

\begin{lem}\label{in}
All elements $\theta_p, p\in M^G$ are linearly independent over
$H^*(BG;{\Bbb Z}_2)$.
\end{lem}
\begin{proof}
Let $\sum_{p\in M^G}l_p\theta_p=0$, where $l_p\in H^*(BG;{\Bbb
Z}_2)$. From \cite[Proposition 5.3.14(2)]{ap}, one knows that
$i^*_q(\theta_p)=0$ for $q\not=p$ in $M^G$, so
$$i^*_q(\sum_{p\in M^G}l_p\theta_p)=\sum_{p\in
M^G}l_pi^*_q(\theta_p)=l_qi^*_q(\theta_q)=l_q\chi_G(q)=0.$$ Since
$\chi_G(q)=t^n$ is a unit in $S^{-1}H^*_G(\{q\};{\Bbb Z}_2)\cong
S^{-1}H^*(BG;{\Bbb Z}_2)$, one has $l_q=0$.
\end{proof}

\begin{lem} \label{l1} Let $\alpha\in S^{-1}H^*_G(M;{\Bbb Z}_2)$. Then
$$\alpha=\sum_{p\in M^G}{{f_p\theta_p}\over{t^n}}$$
where $f_p=S^{-1}i^*_p(\alpha)\in S^{-1}H^*(BG;{\Bbb Z}_2)$.
\end{lem}
\begin{proof}
By \cite[Proposition 5.3.18(1)]{ap}, one has that
$$\alpha=\sum_{p\in M^G}S^{-1}i_{p!}(S^{-1}i^*_p(\alpha)/\chi_G(p)).$$
Since $f_p=S^{-1}i^*_p(\alpha)\in S^{-1}H^*(BG;{\Bbb Z}_2)$, one has
that ${{f_p}\over{\chi_G(p)}}={{f_p}\over{t^n}}\in
S^{-1}H^*(BG;{\Bbb Z}_2)$. Since $S^{-1}i_{p!}$ is a
$S^{-1}H^*(BG;{\Bbb Z}_2)$-algebra homomorphism, one has
$$S^{-1}i_{p!}(S^{-1}i^*_p(\alpha)/\chi_G(p))=S^{-1}i_{p!}({{f_p}\over{t^n}})={{f_p}\over{t^n}}S^{-1}
i_{p!}(1_p)={{f_p}\over{t^n}}i_{p!}(1_p)={{f_p\theta_p}\over{t^n}}$$
so $\alpha=\sum\limits_{p\in M^G}{{f_p\theta_p}\over{t^n}}.$
\end{proof}

\begin{rem} \label{loc}
(i) By  Lemma~\ref{in}, one sees from the formula of Lemma~\ref{l1}
that $\{{{\theta_p}\over {t^n}}|p\in M^G\}$ forms a basis of
$S^{-1}H^*_G(M;{\Bbb Z}_2)$ as a $S^{-1}H^*(BG;{\Bbb Z}_2)$-algebra.

\vskip .2cm
 (ii) In some sense, the formula
$\alpha=\sum\limits_{p\in M^G}{{f_p\theta_p}\over{t^n}}$ explicitly
indicates the
 isomorphism
$S^{-1}i^*:S^{-1}H^*_G(M;{\Bbb Z}_2)\longrightarrow
S^{-1}H^*_G(M^G;{\Bbb Z}_2)$ in Theorem~\ref{local}, which is given
by mapping $\alpha=\sum\limits_{p\in M^G}{{f_p\theta_p}\over{t^n}}$
to $\bigoplus\limits_{p\in M^G}{{f_p}\over{t^n}}$.
\end{rem}

The equivariant Gysin homomorphism of collapsing $M$ to a point
gives the $G$-index of $M$, i.e.,
$$\text{Ind}_G: H^*_G(M;{\Bbb Z}_2)\longrightarrow H^{*-n}(BG;{\Bbb
Z}_2).$$

\begin{thm} \label{formula}
For any $\alpha\in S^{-1}H^*_G(M;{\Bbb Z}_2)$,
$$S^{-1}\text{\rm Ind}_G(\alpha)=\sum_{p\in M^G}{{f_p}\over{t^n}}$$
where $f_p=S^{-1}i^*_p(\alpha)\in S^{-1}H^*(BG;{\Bbb Z}_2)$. In
particular, if $\alpha\in H^*_G(M;{\Bbb Z}_2)$, then
$f_p=i^*_p(\alpha)\in H^*(BG;{\Bbb Z}_2)$ and
\begin{equation} \label{ind}
\text{\rm Ind}_G(\alpha)=\sum_{p\in M^G}{{f_p}\over{t^n}}\in
H^*(BG;{\Bbb Z}_2)={\Bbb Z}_2[t].
\end{equation}
\end{thm}
\begin{proof}
By \cite[Lemma 5.3.19]{ap}, one has that
$\text{Ind}_G(\theta_p)=1_p$, so by Lemma~\ref{l1}
$$S^{-1}\text{Ind}_G(\alpha)=\sum_{p\in
M^G}{{f_pS^{-1}\text{Ind}_G(\theta_p)}\over{t^n}}=\sum_{p\in
M^G}{{f_p\cdot 1_p}\over{t^n}}=\sum_{p\in M^G}{{f_p}\over{t^n}}.$$
The last part of Theorem~\ref{formula} follows immediately since
$H^*(BG;{\Bbb Z}_2)\longrightarrow S^{-1}H^*(BG;{\Bbb Z}_2)$ is
injective.
\end{proof}

\begin{rem} a) It should be pointed out that all arguments in this
section can still be carried out if the action group $G$ is a
2-torus $({\Bbb Z}_2)^k$ of rank $k>1$. In this case, $H^*(BG;{\Bbb
Z}_2)$ is a polynomial algebra ${\Bbb Z}_2[t_1,...,t_k]$ where the
$t_i$'s are one-dimensional generators in $H^1(BG;{\Bbb Z}_2)$, so
that the formula~(\ref{ind}) becomes \begin{equation} \label{ind1}
\text{Ind}_G(\alpha)=\sum_{p\in M^G}{{f_p}\over{\chi_G(p)}}\in {\Bbb
Z}_2[t_1,...,t_k] \end{equation} where $\chi_G(p)$ is a polynomial
of degree $n$ in ${\Bbb Z}_2[t_1,...,t_k]$. Note that related
results can also be found in \cite{d} and \cite{ks}.

b)  The formula (\ref{ind1})  is an analogue of the
Atiyah-Bott-Berlin-Vergne formula for the case $G=T$ (i.e., a
torus), see \cite{ab} and \cite{bv}.
\end{rem}

\section{Equivariant cohomology structure}

In this section, our task is to study the structures of  equivariant
cohomology rings of $G$-manifolds in $\Lambda_n$.

\begin{lem} \label{dim1} Let $M\in \Lambda_n^{2r}$.  Then
$$\dim_{{\Bbb Z}_2}H^i_G(M;{\Bbb Z}_2)=\begin{cases}
\sum\limits_{j=0}^ib_j & \text{ if $i\leq n-1$}\\
2r & \text{ if $i\geq n$.}
\end{cases} $$
\end{lem}
\begin{proof}
Let $$P_s(M_G)=\sum_{i=0}^{\infty}\dim_{{\Bbb Z}_2}H^i_G(M;{\Bbb
Z}_2)s^i$$ be  the equivariant Poincar\'e polynomial of
$H^*_G(X;{\Bbb Z}_2)$.  Since $H^*_G(M;{\Bbb Z}_2)$ is a free
$H^*(BG;{\Bbb Z}_2)$-module, one has that $H^*_G(M;{\Bbb
Z}_2)=H^*(M;{\Bbb Z}_2)\otimes_{{\Bbb Z}_2}H^*(BG;{\Bbb Z}_2)$  so
$$P_s(M_G)=\sum_{i=0}^{\infty}\dim_{{\Bbb Z}_2}H^i_G(M;{\Bbb Z}_2)s^i={1\over {1-s}}\sum_{i=0}^n\dim_{{\Bbb
Z}_2}H^i(M;{\Bbb Z}_2)s^i.$$ Write $b_i=\dim_{{\Bbb Z}_2}H^i(M;{\Bbb
Z}_2)$. Note that $b_i=b_{n-i}$ by Poincar\'e duality and
$b_0=b_n=1$ since $M$ is connected. Then
\begin{align*}
& \quad P_s(M_G)=\sum\limits_{i=0}^{\infty}\dim_{{\Bbb
Z}_2}H^i_G(M;{\Bbb Z}_2)s^i={1\over
{1-s}}\sum\limits_{i=0}^n\dim_{{\Bbb
Z}_2}H^i(M;{\Bbb Z}_2)s^i\\
&= b_0+(b_0+b_1)s+\cdots+
(b_0+b_1+\cdots+b_{n-1})s^{n-1}+(b_0+b_1+\cdots+b_n)(s^n+\cdots)
\\
&=b_0+(b_0+b_1)s+\cdots+
(b_0+b_1+\cdots+b_{n-1})s^{n-1}+2r(s^n+\cdots).
\end{align*}
The lemma then follows from this.
\end{proof}

Let $x=(x_1,...,x_{2r})$ and $y=(y_1,...,y_{2r})$ be two vectors in
$({\Bbb Z}_2)^{2r}$. Define $x\circ y$ by
$$x\circ y=(x_1y_1,...,x_{2r}y_{2r}).$$
Then $({\Bbb Z}_2)^{2r}$ forms a commutative ring with respect to
two operations $+$ and $\circ$. Let
$$\mathcal{V}_{2r}=\big\{x=(x_1,...,x_{2r})^\top\in ({\Bbb
Z}_2)^{2r}\big| |x|=\sum\limits_{i=1}^{2r}x_i=0\big\}.$$ Then it is
easy to see that $\mathcal{V}_{2r}$ is a $(2r-1)$-dimensional
subspace of $({\Bbb Z}_2)^{2r}$, and  there is only  such a subspace
in $({\Bbb Z}_2)^{2r}$. However, generally the operation $\circ$ in
$\mathcal{V}_{2r}$ is obviously not closed.

\vskip .2cm

Given a $M\in \Lambda_n^{2r}$,  one then has that the inclusion $i:
M^G\hookrightarrow M$ induces a monomorphism
$$i^*: H_G^*(M;{\Bbb Z}_2)\longrightarrow ({\Bbb Z}_2)^{2r}[t].$$
By Lemma~\ref{dim1},  there are subspaces $V^M_i$ with $\dim
V^M_i=\sum\limits_{j=0}^ib_j$ $ (i=0,..., n-1)$ of $({\Bbb
Z}_2)^{2r}$ such that
$$H^i_G(M;{\Bbb Z}_2)\cong i^*(H^i_G(M;{\Bbb Z}_2))=\begin{cases}
V_i^Mt^i & \text{ if } i\leq n-1\\
({\Bbb Z}_2)^{2r}t^i & \text{ if } i\geq n \end{cases}$$ where
$V_i^Mt^i=\{vt^i| v\in V_i^M\}$.
\begin{lem} \label{p}
 There are  the following properties:
\begin{enumerate}
\item[(1)] ${\Bbb Z}_2\cong V^M_0\subset V_1^M\subset\cdots\subset V_{n-2}^M\subset V_{n-1}^M=\mathcal{V}_{2r}$,
where $V^M_0$ is generated by $(1,...,1)^\top\in ({\Bbb Z}_2)^{2r}$;
\item[(2)] For $d=\sum\limits_{i=0}^{n-1}i d_i<n$ with
each $d_i\geq0$, $v_{\omega_{d_0}}\circ\cdots\circ
v_{\omega_{d_{n-1}}}\in V^M_{d}$, where
$v_{\omega_{d_i}}=v^{(i)}_1\circ\cdots\circ v^{(i)}_{d_i}$ with each
$v^{(i)}_j\in V^M_i$.
\end{enumerate}
\end{lem}

\begin{proof}
For an element $\alpha\in H_G^*(M;{\Bbb Z}_2)$ of degree $d$, one
has that $i^*(\alpha)=vt^d$ where $v\in ({\Bbb Z}_2)^{2r}$. Since
$i^*=\bigoplus\limits_{p\in M^G}i^*_p$, by Theorem~\ref{formula} one
has that
$$\text{Ind}_G(\alpha)=\sum_{p\in
M^G}{{i^*_p(\alpha)}\over{t^n}}={1\over{t^n}}\sum_{p\in
M^G}i^*_p(\alpha)=|v|t^{d-n}\in {\Bbb Z}_2[t]$$ so if $d<n$, then
$|v|$ must be zero. This means that for each $i<n$,  $V_i^M$ is a
subspace of $\mathcal{V}_{2r}$ and $V_{n-1}^M=\mathcal{V}_{2r}$ is
obvious since $\dim V_{n-1}^M=2r-1$. In particular, when $\alpha=1$
is the identity of $H_G^*(M;{\Bbb Z}_2)$,
$i^*(1)=\bigoplus\limits_{p\in M^G}i_p^*(1)=(1,...,1)^\top\in ({\Bbb
Z}_2)^{2r}$. Thus, $V^M_0\cong {\Bbb Z}_2$ is generated by
$(1,...,1)^\top\in ({\Bbb Z}_2)^{2r}$ since $\dim V_0^M=b_0=1$.
Since $H^*_G(M;{\Bbb Z}_2)=H^*(M;{\Bbb Z}_2)\otimes_{{\Bbb
Z}_2}H^*(BG;{\Bbb Z}_2)$, one has that $(1,...,1)^\top t\in
i^*(H^1_G(M;{\Bbb Z}_2))$. Thus, for any $v\in V_i^M$ with $i<n-1$,
$(vt^i)\circ[(1,...,1)^\top t]=vt^{i+1}\in V_{i+1}^Mt^{i+1}$, so one
has that $v\in V_{i+1}^M$. This completes the proof of
Lemma~\ref{p}(1).

\vskip .2cm

As for the proof of Lemma~\ref{p}(2), for each $v^{(i)}_j\in V^M_i$,
since $i^*: H_G^*(M;{\Bbb Z}_2)\longrightarrow ({\Bbb Z}_2)^{2r}[t]$
is injective, there is a class $\alpha_j^{(i)}$ of degree $i$ in
$H^*_G(M;{\Bbb Z}_2)$ such that $i^*(\alpha_j^{(i)})=v^{(i)}_jt^i$.
Since $i^*=\bigoplus\limits_{p\in M^G}i_p^*$ is also a ring
homomorphism, one has that
$$
i^*(\prod_{i=0}^{n-1}\prod_{j=1}^{d_i}\alpha_j^{(i)})=\bigoplus_{p\in
M^G}i_p^*(\prod_{i=0}^{n-1}\prod_{j=1}^{d_i}\alpha_j^{(i)})
=\bigoplus_{p\in
M^G}\prod_{i=0}^{n-1}\prod_{j=1}^{d_i}i_p^*(\alpha_j^{(i)})=v_{\omega_{d_0}}\circ\cdots\circ
v_{\omega_{d_{n-1}}}t^d
$$
so $v_{\omega_{d_0}}\circ\cdots\circ v_{\omega_{d_{n-1}}}\in V_d^M$.
\end{proof}

\begin{rem} \label{ring} An easy observation shows that  the properties (1) and (2) of Lemma~\ref{p} exactly
give a subring structure of
$$\mathcal{R}_M=V^M_0+V_1^Mt+\cdots+V^M_{n-1}t^{n-1}+({\Bbb Z}_2)^{2r}(t^n+\cdots)$$
in $({\Bbb Z}_2)^{2r}[t]$.
\end{rem}

\vskip .2cm Combining Lemmas~\ref{dim1}, \ref{p}, and
Remark~\ref{ring} one has
\begin{thm} \label{ring str}
Let $M\in \Lambda_n^{2r}$. Then there are subspaces $V^M_i$ with
$\dim V^M_i=\sum\limits_{j=0}^ib_j (i=0,..., n-1)$ of
$\mathcal{V}_{2r}$ such that $H^*_G(M;{\Bbb Z}_2)$ is isomorphic to
the graded ring
$$\mathcal{R}_M=V^M_0+ V^M_1t+\cdots+
V^M_{n-2}t^{n-2}+ V^M_{n-1}t^{n-1}+ ({\Bbb
Z}_2)^{2r}(t^{n}+t^{n+1}+\cdots)$$ where the ring structure of
$\mathcal{R}_M$ is given by
\begin{enumerate}
\item[(a)] ${\Bbb Z}_2\cong V^M_0\subset V_1^M\subset\cdots\subset V_{n-2}^M\subset V_{n-1}^M=\mathcal{V}_{2r}$,
where $V^M_0$ is generated by $(1,...,1)^\top\in ({\Bbb Z}_2)^{2r}$;
\item[(b)] For $d=\sum\limits_{i=0}^{n-1}i d_i<n$ with
each $d_i\geq0$, $v_{\omega_{d_0}}\circ\cdots\circ
v_{\omega_{d_{n-1}}}\in V^M_{d}$, where
$v_{\omega_{d_i}}=v^{(i)}_1\circ\cdots\circ v^{(i)}_{d_i}$ with each
$v^{(i)}_j\in V^M_i$.
\end{enumerate}
\end{thm}

\begin{rem}\label{module}
Since $H^*_G(M;{\Bbb Z}_2)$ is also a free $H^*(BG;{\Bbb
Z}_2)$-module, one has that $$\mathcal{R}_M=V^M_0+ V^M_1t+\cdots+
V^M_{n-2}t^{n-2}+ V^M_{n-1}t^{n-1}+ ({\Bbb
Z}_2)^{2r}(t^{n}+t^{n+1}+\cdots)$$ is a free ${\Bbb Z}_2[t]$-module,
too.
\end{rem}

Next, let us determine the largest-dimensional space $V$ in
 $\{V_i^M| 0\leq i\leq n-1\}$ with the property
that $u\circ v\in \mathcal{V}_{2r}$ for $u,v\in V$.

\vskip .2cm

 Given a vector $v\in \mathcal{V}_{2r}$, set
$$\mathcal{V}(v):=\{x\in \mathcal{V}_{2r}|x\circ v\in
\mathcal{V}_{2r}\}.$$ Then it is easy to see that $\mathcal{V}(v)$
is a linear subspace of $\mathcal{V}_{2r}$, and
$\mathcal{V}(\underline{1}+v)=\mathcal{V}(v)$, and
$\mathcal{V}(v)=\mathcal{V}_{2r}\Longleftrightarrow v=\underline{0}
\text{ or } \underline{1}$ where $\underline{0}=(0,...,0)^\top$ and
$\underline{1}=(1,...,1)^\top$.

\begin{lem} \label{d1}
Let $v\in \mathcal{V}_{2r}$ with $v\neq
\underline{0},\underline{1}$. Then $\dim\mathcal{V}(v) =2r-2$.
\end{lem}
\begin{proof}
Let $\mathcal{V}_1(v)=\{x\in \mathcal{V}_{2r}| v\circ x=x\}$ and
$\mathcal{V}_2(v)=\{x\in \mathcal{V}_{2r} | v\circ x=\underline{0}
\}$. Obviously, they are subspace of $\mathcal{V}_{2r}$ and
$\mathcal{V}_1(v)\cap \mathcal{V}_2(v)=\{\underline{0}\}$. For any
$x\in \mathcal{V}(x)$, $x=v\circ x+(v\circ x+x)$. Since $v\circ x\in
\mathcal{V}_1(v)$ and $v\circ x+x\in \mathcal{V}_2(v)$, one has that
$\mathcal{V}(v)\subset \mathcal{V}_1(v)\oplus \mathcal{V}_2(v)$.
However, obviously $\mathcal{V}_1(v)\oplus \mathcal{V}_2(v)\subset
\mathcal{V}(v)$. Thus,  $\mathcal{V}(v)=\mathcal{V}_1(v)\oplus
\mathcal{V}_2(v)$.

\vskip .2cm

By $\sharp v$ one denotes the number of nonzero elements in $v$. Let
$\sharp v=m$. Then $0<m<2r$ since $v\neq
\underline{0},\underline{1}$, and $m$ is even since $v\in
\mathcal{V}_{2r}$. An easy observation shows that
$\dim\mathcal{V}_1(v)=m-1$ and $\dim\mathcal{V}_2(v)=(2r-m)-1$.
Thus,
$$\dim\mathcal{V}(v)=\dim\mathcal{V}_1(v) +
\dim\mathcal{V}_2(v)=m-1+(2r-m)-1=2r-2.$$
\end{proof}

\begin{lem} \label{d2}
Let $u, v\in \mathcal{V}_{2r}$. Then
$$\mathcal{V}(u+v)=[\mathcal{V}(u)\cap \mathcal{V}(v)]\cup[\mathcal{V}_{2r}\backslash
(\mathcal{V}(u)\cup \mathcal{V}(v))].$$
\end{lem}

\begin{proof} For any $x\in\mathcal{V}(u+v)$, one has that
$x\circ(u+v)=x\circ u+x\circ v\in \mathcal{V}_{2r}$ so $\sharp
(x\circ u+x\circ v)$ is even.  Since $\sharp (a+b)= \sharp a +\sharp
b -2\sharp (a\circ b)$ for any $a, b\in ({\Bbb Z}_2)^{2r}$, one has
that both $\sharp (x\circ u)$ and $\sharp(x\circ v)$ are even or
odd. Thus
\begin{align*}
&\quad \mathcal{V}(u+v)\\&=\{x\in \mathcal{V}_{2r}\vert x\circ u+x\circ v\in \mathcal{V}_{2r}\}\\
&=\{x\in \mathcal{V}_{2r}\vert  \text{ $\sharp (x\circ u)$ and
$\sharp(x\circ v)$ are even}\}\cup\{x\in \mathcal{V}_{2r}\vert  \text{ $\sharp (x\circ u)$ and $\sharp(x\circ v)$ are odd}\}\\
&=[\mathcal{V}(u)\cap \mathcal{V}(v)]\cup[\mathcal{V}_{2r}\backslash
(\mathcal{V}(u)\cup \mathcal{V}(v))].
\end{align*}
\end{proof}

Let $v_1,...,v_k\in \mathcal{V}_{2r}$, and let
$\mathcal{V}(v_1,...,v_k)$ denote $\mathcal{V}(v_1)\cap \cdots\cap
\mathcal{V}(v_k)$.

\begin{prop}\label{dim}
Suppose that $\underline{1},v_1,...,v_k\in \mathcal{V}_{2r}$ are
linearly independent. Then
$$\dim\mathcal{V}(v_1,...,v_k)=2r-1-k.$$
\end{prop}

\begin{proof} One uses induction on $k$. From Lemma~\ref{d1} one knows that  the case
$k=1$ holds. If $k\leq m$, suppose inductively that
Proposition~\ref{dim} holds. Consider the case $k=m+1$. Since
$\mathcal{V}(v_1,...,v_{m+1})\subset \mathcal{V}(v_1,...,v_m)$, one
has $\dim\mathcal{V}(v_1,...,v_{m+1})\leq\dim
\mathcal{V}(v_1,...,v_m)$.

\vskip .2cm If $\dim\mathcal{V}(v_1,...,v_{m+1})=2r-1-m$, then
$\mathcal{V}(v_1,...,v_{m+1})= \mathcal{V}(v_1,...,v_m)$. One claims
that this is impossible. If so, by induction,
$\dim\mathcal{V}(v_2,...,v_{m+1})=2r-1-m$ so
\begin{equation}\label{eq1}
\mathcal{V}(v_2,...,v_{m+1})=\mathcal{V}(v_1,...,v_m) \end{equation}
since
$\mathcal{V}(v_1,...,v_m)=\mathcal{V}(v_1,...,v_{m+1})\subset\mathcal{V}(v_2,...,v_{m+1})$.
In a similar way, one has also that
\begin{equation}\label{eq2}
\mathcal{V}(v_1, v_3,...,v_{m+1})=\mathcal{V}(v_1,...,v_m).
\end{equation}
By Lemma~\ref{d2},
 \begin{align*}\mathcal{V}(v_1+v_2)\cap
\mathcal{V}(v_1)&=\{[\mathcal{V}(v_1)\cap
\mathcal{V}(v_2)]\cup[\mathcal{V}_{2r}\backslash
(\mathcal{V}(v_1)\cup \mathcal{V}(v_2))]\}\cap\mathcal{V}(v_1)\\
&=[\mathcal{V}(v_1)\cap \mathcal{V}(v_2)]\cap \emptyset\\
&=\mathcal{V}(v_1)\cap \mathcal{V}(v_2).
\end{align*}
Similarly, one also has that $\mathcal{V}(v_1+v_2)\cap
\mathcal{V}(v_2)=\mathcal{V}(v_1)\cap \mathcal{V}(v_2)$. Thus
\begin{align*}&\quad \mathcal{V}(v_1+v_2,
v_3,...,v_{m+1})\cap\mathcal{V}(v_1,...,v_m)\\
&=\{\mathcal{V}(v_1+v_2)\cap\mathcal{V}( v_3,...,v_{m+1})\}\cap
\{\mathcal{V}(v_1)\cap
\mathcal{V}(v_2)\cap\mathcal{V}(v_3,...,v_m)\}\\
&= \mathcal{V}(v_1,...,v_m)\end{align*} so
\begin{equation*}
\mathcal{V}(v_1+v_2, v_3,...,v_{m+1})=\mathcal{V}(v_1,...,v_m).
\end{equation*}
On the other hand,
\begin{align*}
&\quad\mathcal{V}(v_1+v_2, v_3,...,v_{m+1})(=\mathcal{V}(v_1,...,v_m)=\mathcal{V}(v_1,...,v_{m+1})) \\
&=\mathcal{V}(v_1+v_2)\cap \mathcal{V}(v_3,...,v_{m+1})\\
&=\{[\mathcal{V}(v_1)\cap
\mathcal{V}(v_2)]\cup[\mathcal{V}_{2r}\backslash
(\mathcal{V}(v_1)\cup \mathcal{V}(v_2))]\}\cap
\mathcal{V}(v_3,...,v_{m+1}) \text{ by Lemma~\ref{d2}}\\
&=\mathcal{V}(v_1,...,v_{m+1})\cup \{[\mathcal{V}_{2r}\backslash
(\mathcal{V}(v_1)\cup \mathcal{V}(v_2))]\cap
\mathcal{V}(v_3,...,v_{m+1})\}.\\
\end{align*}
This implies that
\begin{equation}\label{eq3}
[\mathcal{V}_{2r}\backslash (\mathcal{V}(v_1)\cup
\mathcal{V}(v_2))]\cap \mathcal{V}(v_3,...,v_{m+1})=\emptyset.
\end{equation}
Combining the formulae (\ref{eq1}), (\ref{eq2}) and (\ref{eq3}), one
obtains that
\begin{align*}
\mathcal{V}(v_1,...,v_m)&= [\mathcal{V}(v_1)\cup
\mathcal{V}(v_2)]\cap \mathcal{V}(v_3,...,v_{m+1}) \text{ by
(\ref{eq1}) and (\ref{eq2})}\\
&= \{\mathcal{V}(v_1)\cup
\mathcal{V}(v_2)\cup[\mathcal{V}_{2r}\backslash
(\mathcal{V}(v_1)\cup \mathcal{V}(v_2))]\}\cap
\mathcal{V}(v_3,...,v_{m+1})\text{ by
 (\ref{eq3})}\\
 &= \mathcal{V}_{2r}\cap
\mathcal{V}(v_3,...,v_{m+1})\\
&= \mathcal{V}(v_3,...,v_{m+1})
\end{align*}
so $\mathcal{V}(v_1,...,v_m)=\mathcal{V}(v_3,...,v_{m+1})$. However,
by induction, $\dim\mathcal{V}(v_1,...,v_m)=2r-1-m$ but
$\dim\mathcal{V}(v_3,...,v_{m+1})=2r-1-(m-1)$.  This is a
contradiction.

\vskip .2cm

Therefore,
$\dim\mathcal{V}(v_1,...,v_{m+1})<\dim\mathcal{V}(v_1,...,v_m)$.

\vskip .2cm Now let $\{x_1,..., x_{2r-1-m}\}$ be a basis of
$\mathcal{V}(v_1,...,v_m)$. Since
$$\mathcal{V}(v_1,...,v_{m+1})\subsetneqq \mathcal{V}(v_1,...,v_m)$$
there must exist at least one element $x$ in $\{x_1,...,
x_{2r-1-m}\}$ such that $x\not\in\mathcal{V}(v_1,...,v_{m+1})$. With
no loss, one may assume that
$x_1\not\in\mathcal{V}(v_1,...,v_{m+1})$. Then $x_1\circ
v_{m+1}\not\in \mathcal{V}_{2r}$. If there is also another $x_i
(i\not=1)$ in $\{x_1,x_2,..., x_{2r-1-m}\}$ such that
$x_i\not\in\mathcal{V}(v_1,...,v_{m+1})$, then $x_i\circ
v_{m+1}\not\in \mathcal{V}_{2r}$. However, one knows from the proof
of Lemma~\ref{d2} that $x_1\circ v_{m+1}+x_i\circ v_{m+1}\in
\mathcal{V}_{2r}$ so $x_1+x_i\in \mathcal{V}(v_{m+1})$ and
$x_1+x_i\in \mathcal{V}(v_1,...,v_{m+1})$. In this case, $\{x_1,...,
x_{i-1},x_1+x_i, x_{i+1},..., x_{2r-1-m}\}$ is still a basis of
$\mathcal{V}(v_1,...,v_m)$ but $x_1+x_i\in
\mathcal{V}(v_1,...,v_{m+1})$.
 One can use this way to further modify the basis
 $\{x_1,..., x_{i-1},x_1+x_i, x_{i+1}, ..., x_{2r-1-m}\}$ into a  basis
 $\{x_1, {x'}_2,..., {x'}_{2r-1-m}\}$, such that ${x'}_2,..., {x'}_{2r-1-m}\in
\mathcal{V}(v_1,...,v_{m+1})$ except for $x_1\not\in
\mathcal{V}(v_1,...,v_{m+1})$. Thus,
$\dim\mathcal{V}(v_1,...,v_{m+1})=2r-1-(m+1)$. This completes the
induction and the proof of Proposition~\ref{dim}.
\end{proof}

\begin{cor}\label{d3}
Suppose that $\underline{1},v_1,...,v_k\in \mathcal{V}_{2r}$ are
linearly independent with $v_i\circ v_j\in \mathcal{V}_{2r}$ for any
$i,j\in\{1,...,k\}$. Then $k\leq r-1$.
\end{cor}

\begin{proof}
Since $v_i\circ v_j\in \mathcal{V}_{2r}$ for any
$i,j\in\{1,...,k\}$,
$$\text{Span}\{\underline{1},v_1,...,v_k\}\subset \mathcal{V}(v_1,...,v_k).$$
So $k+1\leq 2r-1-k$, i.e., $k\leq r-1$.
\end{proof}

\begin{cor}\label{largest-dim}
The largest-dimensional space of preserving the operation $\circ$
closed in $\mathcal{V}_{2r}$ of $\{V_i^M| 0\leq i\leq n-1\}$ is
$V^M_{[{n\over 2}]}$ with dimension $r$ if $n$ is odd, and
$V^M_{{n\over 2}-1}$ with dimension $<r$ if $n$ is even but
$b_{{n\over 2}}\not=0$, and $V^M_{{n\over 2}-1}=V^M_{{n\over 2}}$
with dimension $r$ if $n$ is even and $b_{{n\over 2}}=0$.
\end{cor}
\begin{proof}
Since $b_0+b_1+\cdots+b_{n-1}+b_n=2r$ and $b_i=b_{n-i}$ (note
$b_0=b_n=1$), if $n$ is odd, then $b_0+b_1+\cdots+b_{[{n\over
2}]}=r$ so $\dim V^M_{[{n\over 2}]}=r$. Since $\underline{1}\in
V_i^M$ for each $i$, by Corollary~\ref{d3}, the largest-dimensional
space of preserving the operation $\circ$ closed in
$\mathcal{V}_{2r}$ of $\{V_i^M| 0\leq i\leq n-1\}$ must be
$V^M_{[{n\over 2}]}$. If $n$ is even and $b_{{n\over 2}}\not=0$,
then $b_0+b_1+\cdots+b_{{n\over 2}-1}<r$ but
$b_0+b_1+\cdots+b_{{n\over 2}}>r$ so the required
largest-dimensional subspace is $V^M_{{n\over 2}-1}$ with dimension
$r-{{b_{{n\over 2}}}\over 2}<r$. If $n$ is even and $b_{{n\over
2}}=0$, then $$b_0+b_1+\cdots+b_{{n\over
2}-1}=b_0+b_1+\cdots+b_{{n\over 2}}=r$$ so the desired result holds.
\end{proof}

\section{2-dimensional case}

Let $M\in \Lambda_2^{2r}$. By Theorem~\ref{ring str}, we know that
$$H_G^*(M;{\Bbb Z}_2)\cong V^M_0+\mathcal{V}_{2r}t+({\Bbb
Z}_2)^{2r}(t^2+t^3+\cdots).$$ This means that the ring structure of
$H_G^*(M;{\Bbb Z}_2)$ only depends upon the number $2r=|M^G|$. Thus
we have
\begin{prop} \label{2-dim}
Let $M_1, M_2\in \Lambda_2$. Then $H_G^*(M_1;{\Bbb Z}_2)$ and
$H_G^*(M_2;{\Bbb Z}_2)$ are isomorphic if and only if
$|M_1^G|=|M_2^G|$.
\end{prop}

An easy observation shows that  for an orientable connected closed
surface $\Sigma_g$ with genus $g\geq 0$, $\Sigma_g$ must admit a
$G$-action such that $|\Sigma_g^G|=2(g+1)$. Thus, for each $r\geq
1$, $\Lambda_2^{2r}$ is non-empty.

\begin{rem}
However, one knows from \cite{b} that ${\Bbb R}P^2$ never admits a
$G$-action such that the fixed point set is a finite set so ${\Bbb
R}P^2$ doesn't belong to $\Lambda_2$. Actually, generally each
non-orientable connected closed surface $S_g$ with genus $g$ odd
must not belong to $\Lambda_2$. This is because the sum of all mod 2
Betti numbers of $S_g$ is $2+g$, which is odd.
\end{rem}

As a consequence of Proposition~\ref{2-dim}, one has

 \begin{cor}
For each positive integer $r$, all $G$-manifolds in $\Lambda_2^{2r}$
determine a unique equivariant cohomology up to isomorphism.
 \end{cor}

\begin{rem}
For each $M\in \Lambda_2^{2r}$, its equivariant cohomology
$H^*_G(M;{\Bbb Z}_2)$ may be expressed in a simpler way. Since
$(1,...,1)^\top\in V_0^M\cong{\Bbb Z}_2$ and $|v|=0$ for $v\in
\mathcal{V}_{2r}$, one has that
$$H_G^*(M;{\Bbb Z}_2)\cong \Big\{\alpha=(\alpha_1,...,\alpha_{2r})\in({\Bbb
Z}_2)^{2r}[t]\Big| \begin{cases}\alpha_1=\cdots=\alpha_{2r}
&\text{ if $\deg\alpha=0$}\\
\sum\limits_{i=1}^{2r}\alpha_i=0 &\text{ if } \deg\alpha=1
\end{cases}\Big\}.$$ Compare with \cite[Proposition 3.1]{gh}, Goldin
and Holm gave a description for the equivariant cohomolgy of a
compact connected symplectic 4-dimensional manifold with an
effective Hamiltonian $S^1$-action with a finite fixed set, so that
they computed the equivariant cohomology of certain manifolds with a
Hamiltonian action of a torus $T$ (see \cite[Theorem 2]{gh}).
Similarly to the argument in \cite{gh}, using this description of
$H_G^*(M;{\Bbb Z}_2)$ and the main result of Chang and Skjelbred
\cite{cs}, one may give an explicit description in ${\Bbb
Z}_2[t_1,...,t_k]$ of the equivariant cohomology of a closed $({\Bbb
Z}_2)^k$-manifold $N$ with the following conditions:
\begin{enumerate}
 \item[(1)] The fixed point set is finite;
 \item[(2)] The equivariant cohomology of $N$ is a free $H^*(B({\Bbb
 Z}_2)^k;{\Bbb Z}_2)$-module;
 \item[(3)] For $K< ({\Bbb Z}_2)^k$ a corank-1 2-torus,  each component
 of $N^K$ has dimension at most 2.
\end{enumerate}
We would like to leave it to readers as an exercise. For this
description in ${\Bbb Z}_2[t_1,...,t_k]$ of $H^*_{({\Bbb
Z}_2)^k}(N;{\Bbb Z}_2)$, actually the above restriction condition
(3) for each component
 of $N^K$ is the best possible since generally
there can be different equivariant cohomology structures for
$G$-manifolds in $\Lambda_n^{2r}$ when $n\geq 3$ (see, e.g., Section
6 of this paper).
\end{rem}

\section{An analytic description of ring isomorphisms and a necessary and sufficient condition}

The purpose of this section is to give an analytic description for
the isomorphism between equivariant cohomology rings of two
$G$-manifolds in $\Lambda_n$ and to show  a necessary and sufficient
condition that the equivariant cohomology rings of two $G$-manifolds
in $\Lambda_n$ are isomorphic.

\vskip .2cm

 Let
$$\mathcal{W}_{2r}=\{\sigma\in {\a} \mathcal{V}_{2r}\subset \G (2r,{\Bbb Z}_2)|
\text{ $\sigma (x\circ y)=\sigma (x)\circ\sigma (y)$ for any  $x,
y\in \mathcal{V}_{2r}$}\}.$$ Given two $\sigma,\tau$ in
$\mathcal{W}_{2r}$, it is easy to check that $\sigma\tau\in
\mathcal{W}_{2r}$. Thus one has that
 $\mathcal{W}_{2r}$ is a subgroup of $\a \mathcal{V}_{2r}$.

\begin{lem}\label{weyl}
$\mathcal{W}_{2r}$ is the Weyl subgroup of $\G(2r,{\Bbb Z}_2)$.
\end{lem}

\begin{proof}
By \cite{ab1}, it suffices to prove that $\mathcal{W}_{2r}$ is
isomorphic to the symmetric group $\mathcal{S}_{2r}$ of rank $2r$.
Let $\sigma=(a_{ij})_{2r\times 2r}$ be an element of
$\mathcal{W}_{2r}\subset \G(2r,{\Bbb Z}_2)$. Since $\sigma$ is an
automorphism of $\mathcal{V}_{2r}$, there exists a vector
$\underline{x}\in\mathcal{V}_{2r}$ such that
$\sigma(\underline{x})=\underline{1}$ where
$\underline{1}=(1\underbrace{,...,}_{2r}1)^\top$. Since
$\underline{1}\circ x=x$ for any $x\in\mathcal{V}_{2r}$, one has
that
$\underline{1}=\sigma(\underline{x})=\sigma(\underline{1}\circ\underline{x})=\sigma(\underline{1})\circ
\sigma(\underline{x})=\sigma(\underline{1})$ so
\begin{equation} \label{e1}
\sum_{j=1}^{2r}a_{ij}=1, \text{ for } i=1,2,...,2r.
\end{equation} On the other
hand, for any $x=(x_1,...,x_{2r})^\top$ and
$y=(y_1,...,y_{2r})^\top$ in $\mathcal{V}_{2r}$, one has that
$\sigma(x\circ y)=\sigma(x)\circ\sigma(y)$ so
\begin{equation}\label{e2}
\sum_{l=1}^{2r}a_{il}x_ly_l=(\sum_{j=1}^{2r}a_{ij}x_j)(\sum_{k=1}^{2r}a_{ik}y_k),
\text{ for } i=1,2,...,2r.
\end{equation}
From (\ref{e1}) one knows that for each $i$, the number $q(i)$ of
nonzero elements in $a_{i1},...,a_{i2r}$ is odd.

\vskip .2cm

Now let us show that for each $i$, $q(i)$ actually must be 1. Taking
an $i$, without loss of generality one may assume that
$a_{i1}=\cdots=a_{iq(i)}=1$, and $a_{i(q(i)+1)}=\cdots=a_{i2r}=0$.
Then from (\ref{e2}) one has
\begin{equation}\label{e3}
(\sum_{j=1}^{q(i)}x_j)(\sum_{k=1}^{q(i)}y_k)+\sum_{l=1}^{q(i)}x_ly_l=0.
\end{equation}
If $q(i)>1$, taking $x$ with $x_1=x_{q(i)}=1$ and $x_j=0$ for
$j\not=1, q(i)$ and $y$ with $y_2=y_{q(i)}=1$ and $y_k=0$ for
$k\not=2,q(i)$, the left side of (\ref{e3}) then becomes 1, but this
is impossible. Thus $q(i)=1$.

\vskip .2cm

Since $q(i)=1$ for each $i$, this means that $\sigma$ is actually
obtained by doing a permutation on all rows (or all columns) of the
identity matrix. The lemma then follows from this.
\end{proof}

\begin{thm}\label{analyze}
Let $M_1$ and $M_2$ in $\Lambda_n^{2r}$. Suppose that $f$ is an
isomorphism  between graded rings
$$\mathcal{R}_{M_1}=V^{M_1}_0+\cdots+ V^{M_1}_{n-2}t^{n-2}+\mathcal{V}_{2r}t^{n-1}+ ({\Bbb
Z}_2)^{2r}(t^n+\cdots)$$ and
$$\mathcal{R}_{M_2}=V^{M_2}_0+\cdots+ V^{M_2}_{n-2}t^{n-2}+\mathcal{V}_{2r}t^{n-1}+ ({\Bbb
Z}_2)^{2r}(t^n+\cdots).$$  Then there is an element $\sigma\in
\mathcal{W}_{2r}$ such that $f=\sum\limits_{i=0}^\infty \sigma t^i$
is analytic, where  $f(\beta)=\sum\limits_{i=0}^\infty \sigma(v_i)
t^i$ for $\beta=\sum\limits_{i=0}^\infty v_i t^i\in
\mathcal{R}_{M_1}$.
\end{thm}

\begin{proof}
Since the restriction
$f|_{\mathcal{V}_{2r}t^{n-1}}:\mathcal{V}_{2r}t^{n-1}\longrightarrow\mathcal{V}_{2r}t^{n-1}$
is a linear isomorphism, there exists an automorphism $\sigma$ of $
\mathcal{V}_{2r}$ such that $f|_{\mathcal{V}_{2r}t^{n-1}}=\sigma
t^{n-1}$. \vskip .2cm

First, let us show that $\sigma\in \mathcal{W}_{2r}$.  Since $x\circ
x=x$ for any $x\in ({\Bbb Z}_2)^{2r}$ and  $f$ is a ring
isomorphism, one has that for any $v\in \mathcal{V}_{2r}$,
$$f(vt^{2n-2})=f(v\circ
vt^{2n-2})=f(vt^{n-1})f(vt^{n-1})=\sigma(v)\circ\sigma(v)t^{2n-2}=\sigma(v)t^{2n-2}$$
so $f|_{\mathcal{V}_{2r}t^{2n-2}}=\sigma t^{2n-2}$.  Furthermore,
for any $v_1, v_2\in \mathcal{V}_{2r}$, one has that
$$\sigma(v_1\circ v_2)t^{2n-2}=f(v_1\circ
v_2t^{2n-2})=f(v_1t^{n-1})f(v_2t^{n-1})=\sigma(v_1)\circ\sigma(v_2)t^{2n-2}$$
so $\sigma(v_1\circ v_2)=\sigma(v_1)\circ\sigma(v_2)$. Thus,
$\sigma\in \mathcal{W}_{2r}$.

\vskip .2cm

By Theorem~\ref{ring str}, for each $i<n-1$, $\underline{1}\in
V_i^M\subset \mathcal{V}_{2r}$, and  one knows from the proof of
Lemma~\ref{weyl} that $\sigma(\underline{1})=\underline{1}$ so
$f|_{V_0^M}=\sigma$. By Remark~\ref{module} one knows that
$\mathcal{R}_{M_i}, i=1,2,$ are free ${\Bbb Z}_2[t]$-modules, and it
is easy to see that $f$ is also an isomorphism between  free ${\Bbb
Z}_2[t]$-modules $\mathcal{R}_{M_1}$ and $\mathcal{R}_{M_2}$. For
$0\leq i<n-1$, let $v\in V_i^M\subset \mathcal{V}_{2r}$, since
$f|_{\mathcal{V}_{2r}t^{n-1}}=\sigma t^{n-1}$, one has that
$$f(vt^i)t^{n-1-i}=f(vt^{n-1})=
\sigma(v)t^{n-1}$$ so $f(vt^i)=\sigma(v)t^i$. Thus,
$f|_{V_i^Mt^i}=\sigma t^i$ for $0\leq i<n-1$.

\vskip .2cm

Let $\ell=a(n-1)+b\geq n$ with $b<n-1$. For any $v\in
\mathcal{V}_{2r}$, one has that
$$f(vt^\ell)=f(v\underbrace{\circ \cdots\circ}_avt^{a(n-1)})t^b
=[f(vt^{n-1})]^at^b=\sigma(v)\underbrace{\circ
\cdots\circ}_a\sigma(v)t^\ell=\sigma(v)t^\ell.$$ Thus, for $i\geq
n$, $f|_{\mathcal{V}_{2r}^Mt^i}=\sigma t^i$.

\vskip .2cm  Since $\mathcal{V}_{2r}$ is not closed with respect to
the operation $\circ$ by Corollary~\ref{d3}, there must be $u,
w\in\mathcal{V}_{2r}$ such that $u\circ w\not\in \mathcal{V}_{2r}$.
Since $\mathcal{V}_{2r}$ has dimension $2r-1$,  one then has that
$({\Bbb Z}_2)^{2r}=\mathcal{V}_{2r}+\text{Span}\{u\circ w\}$.
Actually, this is a direct sum decomposition of $({\Bbb Z}_2)^{2r}$,
i.e., $({\Bbb Z}_2)^{2r}=\mathcal{V}_{2r}\oplus\text{Span}\{u\circ
w\}$. Furthermore, let $x\in({\Bbb Z}_2)^{2r}$, then there is a
vector $v\in \mathcal{V}_{2r}$ such that $x$ can be written as
$v+\varepsilon u\circ w$  where  $\varepsilon=0$ or 1. For
$i=a(n-1)+b\geq n$ with $b<n-1$, since
$f|_{\mathcal{V}_{2r}^Mt^j}=\sigma t^j$ for any $j\geq n-1$, one has
that
\begin{align*}
f(xt^i)t^{n-1-b}&=f( vt^i+\varepsilon u\circ wt^i)t^{n-1-b}\\&=
f(vt^i)t^{n-1-b}
+\varepsilon f(u\circ wt^i)t^{n-1-b}\\
&= f(vt^{(a+1)(n-1)})+\varepsilon f([ut^{n-1}]\circ [wt^{a(n-1)}])\\
&= \sigma(v)t^{(a+1)(n-1)}+\varepsilon f(ut^{n-1})f(wt^{a(n-1)})\\
&=\sigma(v)t^{(a+1)(n-1)}+\varepsilon \sigma(u)\circ
\sigma(w)t^{(a+1)(n-1)}\\
&= [\sigma(v)+\varepsilon \sigma(u)\circ
\sigma(w)]t^{(a+1)(n-1)}\\
&= \sigma(v+\varepsilon u\circ
w)t^{(a+1)(n-1)}\text{ since $\sigma\in\mathcal{W}_{2r}$}\\
&= \sigma(x)t^{i+(n-1-b)}\\
\end{align*}
 so $f(xt^i)=\sigma(x)t^i$. Thus
$f|_{({\Bbb Z}_2)^{2r}t^i}=\sigma t^i$ for $i\geq n$.

\vskip .2cm Combining the above argument, we complete the proof.
\end{proof}

\begin{thm} \label{ns}
Let $M_1$ and $M_2$ in $\Lambda_n^{2r}$. Then $H^*_G(M_1;{\Bbb
Z}_2)$ and $H^*_G(M_2;{\Bbb Z}_2)$ are isomorphic if and only if
there exists an element $\sigma\in \mathcal{W}_{2r}$ such that
$\sigma$ isomorphically  maps $V^{M_1}_i$ onto $V^{M_2}_i$ for $i<
n-1$.
\end{thm}

\begin{proof}
Suppose that $H^*_G(M_1;{\Bbb Z}_2)$ and $H^*_G(M_2;{\Bbb Z}_2)$ are
isomorphic. Then by Theorem~\ref{ring str} there is an isomorphism
$f$ between graded rings
$$\mathcal{R}_{M_1}=V^{M_1}_0+\cdots+ V^{M_1}_{n-2}t^{n-2}+\mathcal{V}_{2r}t^{n-1}+ ({\Bbb
Z}_2)^{2r}(t^n+\cdots)$$ and
$$\mathcal{R}_{M_2}=V^{M_2}_0+\cdots+ V^{M_2}_{n-2}t^{n-2}+\mathcal{V}_{2r}t^{n-1}+ ({\Bbb
Z}_2)^{2r}(t^n+\cdots).$$  One knows from Theorem~\ref{analyze} that
there is an element $\sigma\in \mathcal{W}_{2r}$ such that
$f=\sum\limits_{i=0}^\infty \sigma t^i$. Then the restriction
$f|_{V_i^{M_1}t^i}=\sigma t^i$, which  is an isomorphism from
$V^{M_1}_it^i$ to $V^{M_2}_it^i$ for $i< n-1$. Thus $\sigma$
isomorphically maps $V^{M_1}_i$ onto $V^{M_2}_i$ for $i< n-1$.

\vskip .2cm

Conversely, if   there exists an element $\sigma\in
\mathcal{W}_{2r}$ such that $\sigma$ isomorphically  maps
$V^{M_1}_i$ onto $V^{M_2}_i$ for $i< n-1$, then by
Theorem~\ref{analyze}, $\sum\limits_{i=0}^\infty \sigma t^i$ gives
an isomorphism between graded rings
$\mathcal{R}_{M_1}=V^{M_1}_0+\cdots+
V^{M_1}_{n-2}t^{n-2}+\mathcal{V}_{2r}t^{n-1}+ ({\Bbb
Z}_2)^{2r}(t^n+\cdots)$ and $\mathcal{R}_{M_2}=V^{M_2}_0+\cdots+
V^{M_2}_{n-2}t^{n-2}+\mathcal{V}_{2r}t^{n-1}+ ({\Bbb
Z}_2)^{2r}(t^n+\cdots)$. Then $H^*_G(M_1;{\Bbb Z}_2)$ and
$H^*_G(M_2;{\Bbb Z}_2)$ are isomorphic by Theorem~\ref{ring str}.
\end{proof}

\begin{rem}
Let $M_1, M_2\in \Lambda_n$. One sees from Theorem~\ref{ns} that if
$|M_1^G|\not=|M_2^G|$, then $H_G^*(M_1;{\Bbb Z}_2)$ and
$H_G^*(M_2;{\Bbb Z}_2)$ must not be isomorphic.
\end{rem}

In the  case $n=3$, for $M\in \Lambda_3^{2r}$, since
$$H_G^*(M;{\Bbb Z}_2)\cong V_0^M+V_1^Mt+\mathcal{V}_{2r}t^2+({\Bbb
Z}_2)^{2r}(t^3+\cdots),$$ one sees that the structure of
$H_G^*(M;{\Bbb Z}_2)$  actually depends upon that of $V_1^M$.  Thus,
Theorem~\ref{ns} has a simpler expression in this case.

\begin{cor}\label{3ns}
Let $M_1$ and $M_2$ in $\Lambda_3^{2r}$. Then $H^*_G(M_1;{\Bbb
Z}_2)$ and $H^*_G(M_2;{\Bbb Z}_2)$ are isomorphic if and only if
there exists an element $\sigma\in \mathcal{W}_{2r}$ such that
$\sigma$ isomorphically  maps $V^{M_1}_1$ onto $V^{M_2}_1$.
\end{cor}

\section{The number of equivariant cohomology structures}

In this section we shall consider the number of equivariant
cohomology rings up to isomorphism of all 3-dimensional
$G$-manifolds in $\Lambda_3^{2r}$. For $M\in \Lambda_3^{2r}$ one
knows  that $ V_1^M$ has dimension $r$ and is the
largest-dimensional subspace of $\mathcal{V}_{2r}$ with the property
that $u\circ v\in \mathcal{V}_{2r}$ for $u, v\in V_1^M$ by
Corollary~\ref{largest-dim}.

\vskip .2cm

Let $\mathcal{M}_r$ denote the set of those $2r\times r$ matrices
$(v_1,...,v_r)$ with rank $r$ over ${\Bbb Z}_2$ such that $v_i\circ
v_j\in \mathcal{V}_{2r}$ for any $1\leq i,j\leq r$. $\mathcal{M}_r$
admits the following two actions.

\vskip .2cm

One action is the right action of $\text{GL}(r,{\Bbb Z}_2)$ on
$\mathcal{M}_r$ defined by $(v_1,...,v_r)\lambda$ for
$\lambda\in\text{GL}(r,{\Bbb Z}_2)$. It is easy to see that such
action is free. Obviously, all column vectors of each matrix
$(v_1,...,v_r)$ span the same linear space as all column vectors of
$(v_1,...,v_r)\lambda$ for $\lambda\in\text{GL}(r,{\Bbb Z}_2)$.

\vskip .2cm

The other action is the left action of the Weyl group
$\mathcal{W}_{2r}=\mathcal{S}_{2r}$ on $\mathcal{M}_r$ defined by
$\sigma(v_1,...,v_r)=(\sigma v_1,...,\sigma v_r)$ for
$(v_1,...,v_r)\in \mathcal{M}_r$ and $\sigma\in \mathcal{S}_{2r}$.
In general, this action is not free.

\vskip .2cm

For  $M\in \Lambda_3^{2r}$, since the space
$V_1^M\subset\mathcal{V}_{2r}$  may be spanned by all column vectors
of some matrix $(v_1,...,v_r)$ in $\mathcal{M}_r$,  the number of
all possible spaces $V_1^M\subset\mathcal{V}_{2r}$ is at most
$\big|\mathcal{M}_r/\G(r,{\Bbb
Z}_2)\big|={{|\mathcal{M}_r|}\over{|\G(r,{\Bbb Z}_2)|}}$.

\vskip .2cm

Together with the above understood and Corollary~\ref{3ns}, one has
 \begin{prop}\label{number}
The  number of equivariant cohomology rings  up to isomorphism of
all ${\Bbb Z}_2$-manifolds in $\Lambda_3^{2r}$ is at most
$$|\mathcal{S}_{2r}\backslash\mathcal{M}_r/\G(r,{\Bbb Z}_2)|.$$
 \end{prop}

\begin{rem} \label{vp}
There is a natural map $g$ from $\Lambda_3^{2r}$ to
$\mathcal{M}_r/\text{GL}(r,{\Bbb Z}_2)$. If this map is surjective,
then the number of equivariant cohomology rings  up to isomorphism
of all ${\Bbb Z}_2$-manifolds in $\Lambda_3^{2r}$ is exactly
$|\mathcal{S}_{2r}\backslash\mathcal{M}_r/\text{GL}(r,{\Bbb Z}_2)|.$
To determine whether $g$ is surjective or not is an interesting
thing, but it seems to be quite difficult.
\end{rem}

\vskip .2cm

Generally, the computation of the number
$|\mathcal{S}_{2r}\backslash\mathcal{M}_r/\text{GL}(r,{\Bbb Z}_2)|$
is not an easy thing. Next, we shall analyze  this number.

\vskip .2cm

 Taking an element $A\in\mathcal{M}_r$, there must be
$\sigma\in \mathcal{S}_{2r}$ and $\lambda\in\text{GL}(r,{\Bbb Z}_2)$
such that
$$\sigma A\lambda=
\begin{pmatrix}
I_r \\
P
 \end{pmatrix}$$ so each orbit of the orbit set $\mathcal{S}_{2r}\backslash\mathcal{M}_r/\text{GL}(r,{\Bbb
 Z}_2)$ contains the representative of the form $
\begin{pmatrix}
I_r \\
P
 \end{pmatrix}$, where $I_r$ is the $r\times r$ identity matrix.
 It is easy to check that $P\in O(r,{\Bbb Z}_2)$, where $O(r,{\Bbb Z}_2)$ is the orthogonal matrix group over
 ${\Bbb Z}_2$. Obviously, $O(r,{\Bbb Z}_2)$ always admits the left and
 right actions of the Weyl group $\mathcal{S}_r$ in $\text{GL}(r,{\Bbb
 Z}_2)$.

\begin{lem}  \label{bound}
$|\mathcal{S}_{2r}\backslash\mathcal{M}_r/\G(r,{\Bbb Z}_2)|\leq
|\mathcal{S}_r\backslash O(r,{\Bbb Z}_2)/ \mathcal{S}_r|.$
\end{lem}

 \begin{proof}
Suppose that $P_1,P_2\in O(r,{\Bbb Z}_2)$ belong to the same orbit
in $\mathcal{S}_r\backslash O(r,{\Bbb Z}_2)/ \mathcal{S}_r$. Then
there exist $\tau, \rho$ in $\mathcal{S}_r$ such that
$$\tau P_1\rho=P_2.$$
Let $\sigma=\begin{pmatrix}
\rho^{-1} & 0 \\
0 & \tau
 \end{pmatrix}$. Then $\sigma\in \mathcal{S}_{2r}$. Obviously,
 $$\sigma\begin{pmatrix}
I_r \\
P_1
 \end{pmatrix}\rho=\begin{pmatrix}
\rho^{-1}\rho \\
\tau P_1\rho
 \end{pmatrix}=\begin{pmatrix}
I_r \\
P_2
 \end{pmatrix}$$
 so $
\begin{pmatrix}
I_r \\
P_1
 \end{pmatrix}$ and $
\begin{pmatrix}
I_r \\
P_2
 \end{pmatrix}$ belong to the same orbit in $\mathcal{S}_{2r}\backslash\mathcal{M}_r/\text{GL}(r,{\Bbb
 Z}_2)$. This means that
$$|\mathcal{S}_{2r}\backslash\mathcal{M}_r/\text{GL}(r,{\Bbb Z}_2)|\leq |\mathcal{S}_r\backslash
O(r,{\Bbb Z}_2)/ \mathcal{S}_r|.$$ This completes the proof.
 \end{proof}

\begin{rem}
Generally, $|\mathcal{S}_{2r}\backslash\mathcal{M}_r
/\text{GL}(r,{\Bbb Z}_2)|$ is not equal to $|\mathcal{S}_r\backslash
O(r,{\Bbb Z}_2)/ \mathcal{S}_r|.$ For example, take
$$P_1=\begin{pmatrix}
0&1&1&1&1&1\\1&0&1&1&1&1\\1&1&0&1&1&1\\1&1&1&0&1&1\\1&1&1&1&0&1\\1&1&1&1&1&0
\end{pmatrix} \text{ and } P_2=\begin{pmatrix}
0&1&1&1&1&1\\1&1&0&0&0&1\\1&0&1&0&0&1\\1&0&0&1&0&1\\1&0&0&0&1&1\\1&1&1&1&1&0
\end{pmatrix}$$
it is easy to check that $P_1$ and $P_2$ don't belong to the same
orbit in $\mathcal{S}_6\backslash O(6,{\Bbb Z}_2)/ \mathcal{S}_6$,
but $\begin{pmatrix}
I_6 \\
P_1
 \end{pmatrix}$ and $\begin{pmatrix}
I_6 \\
P_2
 \end{pmatrix}$ belong to the same orbit in
 $\mathcal{S}_{12}\backslash\mathcal{M}_6
/\text{GL}(6,{\Bbb Z}_2)$.
\end{rem}

\begin{defn}
Let $A=(v_1,...,v_r)\in \mathcal{M}_r$. One says that $A$ is {\em
irreducible} if the space $\text{Span}\{v_1,...,v_r\}$ cannot be
decomposed as a direct sum of some nonzero subspaces $V_1,...,V_l,
l>1,$ with the property that $V_i\circ V_j=\{\underline{0}\}$ for
$i\not=j$, where $V_i\circ V_j=\{\underline{0}\}$ means that for any
$x\in V_i$ and $y\in V_j$, $x\circ y=\underline{0}$.
\end{defn}

We would like to point out that if one can  find out all possible
irreducible matrices of $\mathcal{M}_r$ for any $r$, then one  can
construct a representative of each class in
$\mathcal{S}_{2r}\backslash\mathcal{M}_r/\text{GL}(r,{\Bbb Z}_2)$,
with the form $$
\begin{pmatrix}
A_1&0&...&0\\0&A_2&...&0\\
& &...
\\0&0&...&A_s
\end{pmatrix} $$
where the  blocks $A_i$'s are irreducible matrices.

\vskip .2cm

The following  result shows that there exists at least one
irreducible matrix in $\mathcal{M}_r$ for almost any positive
integer $r$. Let
$$A(l)=\left(
\begin{matrix}
1&0&0&0&......&0&0\\
1&0&0&0&......&0&1\\
1&0&0&0&......&1&0\\
&&&&......
\\1&0&0&1&......&0&0\\
1&0&1&0&......&0&0\\
1&1&0&0&......&0&0
\\1&1&1&1&......&1&1\\
0&1&1&1&......&1&1
\\0&1&0&0&......&0&0\\
0&0&1&0&......&0&0\\
0&0&0&1&......&0&0\\
&&&&......
\\0&0&0&0&......&1&0\\
0&0&0&0&......&0&1
\end{matrix} \right)_{(2l-1)\times (l-1)}$$
be a $(2l-1)\times (l-1)$ matrix with $l\geq 4$ and
$\underline{1}_{l\times 1}=(1\underbrace{,...,}_l1)^\top$. Note that
only when $l\geq 4$ is even,
 the first column vector $v$ of $A(l)$ exactly has the property $|v|=0$ in
 ${\Bbb Z}_2$.
\begin{lem} \label{ir}
$(a)$ For  even $r\geq 4$, there exists an irreducible $2r\times r$
matrix $$\begin{pmatrix} A(r)&\underline{1}_{(2r-1)\times 1}\\0&1
\end{pmatrix}.$$ $(b)$ For  odd  $r\geq 7$, there exist irreducible
$2r\times r$ matrices of the following form
$$
\begin{pmatrix}A(s)&0&\underline{1}_{(2s-1)\times 1}\\0&A(t)&\underline{1}_{(2t-1)\times 1}
\end{pmatrix} $$
with only two blocks $A(l)$'s,  where $s+t=r+1$ and $s,t\geq4$ are
even. In particular, both $
\begin{pmatrix}A(s_1)&0&\underline{1}_{(2s_1-1)\times1}\\0&A(t_1)&\underline{1}_{(2t_1-1)\times 1}
\end{pmatrix}$ and $
\begin{pmatrix}A(s_2)&0&\underline{1}_{(2s_2-1)\times1}\\0&A(t_2)&\underline{1}_{(2t_2-1)\times1}
\end{pmatrix}$ belong to the same orbit in $\mathcal{S}_{2r}\backslash\mathcal{M}_r/\G(r,{\Bbb Z}_2)$
if and only if $\{s_1,t_1 \} = \{s_2,t_2 \}$, where $s_j+t_j=r+1$
and $s_j,t_j\geq 4$ are even for $ j=1,2$.
\end{lem}

\begin{rem}
(i) A direct observation shows that when $r=1$, $\mathcal{M}_1$
contains a unique matrix $\begin{pmatrix} 1\\1
\end{pmatrix}$, which is irreducible, and when $r=2,3,5$, there is no any irreducible matrix.
However,  we don't know whether those irreducible matrices stated in
Lemma~\ref{ir}, with $\begin{pmatrix} 1\\1
\end{pmatrix}$ together, give all possible irreducible matrices.

\vskip .2cm   (ii) For odd $r\geq7$, let $\lambda(r)$ denote the
number of the orbit classes of irreducible matrices in
$\mathcal{S}_{2r}\backslash\mathcal{M}_r/\text{GL}(r,{\Bbb Z}_2)$.
By Lemma~\ref{ir}, $\lambda(r)$ is equal to or more than the number
of the solutions $(x, y)$ of the following equation
$$\begin{cases}
\ x+y={{r+1}\over 2}\\
\ x\geq y\geq 2.
\end{cases}$$
Then, an easy argument shows that the number of solutions is
$[{{{(r+1)/2}-1-1}\over2}]=[{{r-3}\over4}]$, so  $\lambda(r)\geq
[{{r-3}\over4}]$.
\end{rem}

Let $A=(v_1,...,v_r)\in \mathcal{M}_r$. Set
$$\mathcal{V}(A):=\text{Span}\{v_1,...,v_r\}$$
and
$$\mathcal{X}(A):=\text{Span}\{v_i\circ v_j| i,j=1,...,r\}.$$
It is easy to see the following properties:
\begin{enumerate}
 \item[(P1)] $\mathcal{V}(A)\subset\mathcal{X}(A)\subset
 \mathcal{V}_{2r}$ so $r\leq \dim\mathcal{X}(A)\leq 2r-1$;
 \item[(P2)] If $A, B\in \mathcal{M}_r$ belong to the same orbit in
  $\mathcal{S}_{2r}\backslash\mathcal{M}_r/\text{GL}(r,{\Bbb Z}_2)$,
  then $\mathcal{X}(A)$ is linearly isomorphic to $\mathcal{X}(B)$.
\end{enumerate}

\vskip .2cm

\noindent {\bf Fact 1.} {\em Let $A\in \mathcal{M}_r$.  If $A$ is
not irreducible, then there exists a $2k\times k$ matrix
$A_{2k\times k}$ and a $2l\times l$ matrix $A_{2l\times l}$ with
$k+l=r$ such that $A$ and $\begin{pmatrix} A_{2k\times k}&0\\0&
A_{2l\times l}
\end{pmatrix}$ belong to the same orbit in
$\mathcal{S}_{2r}\backslash\mathcal{M}_r/\G(r,{\Bbb Z}_2)$.}
\begin{proof}
If $A$ is not irreducible, by definition, there exists a $k_1\times
k_2$ matrix $A_{k_1\times k_2}$ and a $l_1\times l_2$ matrix
$A_{l_1\times l_2}$ with $k_1+l_1=2r$ and $k_2+l_2=r$ such that $A$
and $\begin{pmatrix} A_{k_1\times k_2}&0\\0& A_{l_1\times l_2}
\end{pmatrix}$ belong to the same orbit in
$\mathcal{S}_{2r}\backslash\mathcal{M}_r/\text{GL}(r,{\Bbb Z}_2)$.
By Proposition~\ref{dim} and the proof method of Corollary~\ref{d3},
one has that $k_2\leq k_1/2$ and $l_2\leq l_1/2$. Then the relations
$k_1+l_1=2r$ and $k_2+l_2=r$ force $k_2= k_1/2$ and $l_2=l_1/2$.
\end{proof}

\noindent {\bf Fact 2.} {\em  Let $A\in \mathcal{M}_r$.  If
$\dim\mathcal{X}(A)= 2r-1$, then $A$ is irreducible.}

\begin{proof}
Suppose that  $A$ is not irreducible. By Fact 1, there is a matrix
$\begin{pmatrix} A_{2k\times k}&0\\0& A_{2l\times l}
\end{pmatrix}$ that is in the same orbit as $A$.
Thus $$\dim \mathcal{X}(A)=\dim \mathcal{X}(A_{2k\times k})+\dim
\mathcal{X}(A_{2l\times l})
 \leq (2k-1)+(2l-1)=2r-2$$ by (P1). This is a contradiction.
\end{proof}

\begin{rem}
 Let $A\in \mathcal{M}_r$. Then there is a matrix $P\in O(r, {\Bbb Z}_2)$
such that $A$ and $\begin{pmatrix}I_r\\P
\end{pmatrix}$ belong to  the same orbit in
$\mathcal{S}_{2r}\backslash\mathcal{M}_r/\G(r,{\Bbb Z}_2)$.
Furthermore, an easy argument shows  by (P2) that if
$\dim\mathcal{X}(A)= r$,
then  $P$ can be chosen as being $I_r$, 
so $A$
is not irreducible when $r\not=1$.
\end{rem}

Now let us give the proof of Lemma~\ref{ir}.

\vskip .2cm

{\em Proof of Lemma~\ref{ir}}.  Let $A=\begin{pmatrix}
A(r)&\underline{1}_{(2r-1)\times 1}\\0&1
\end{pmatrix}=(v_1,...,v_{r-1}, \underline{1}_{2r\times1})$.  By Fact 2, it suffices to check
that $\dim\mathcal{X}(A)=2r-1$. A direct observation shows that
$$\{v_1\circ v_2, v_1\circ v_3,...,v_1\circ v_{r-1},v_2\circ v_3,
v_1,v_2,v_3,...,v_{r-1},\underline{1}_{2r\times1} \}$$ are linearly
independent, so $\dim\mathcal{X}(A)\geq 2r-1$. Since
$\dim\mathcal{X}(A)\leq 2r-1$ by (P1), $\dim\mathcal{X}(A)=2r-1$.
This completes the  proof of Lemma~\ref{ir}(a).

\vskip .2cm

 In a similar way to the above argument,  one may show that $
\begin{pmatrix}A(s)&0&\underline{1}_{(2s-1)\times 1}\\0&A(t)&\underline{1}_{(2t-1)\times 1}
\end{pmatrix}$ is irreducible, too.

\vskip .2cm

Suppose that $
\begin{pmatrix}A(s_1)&0&0&...&0&\underline{1}_{(2s_1-1)\times 1}\\0&A(s_2)&0&...&0&\underline{1}_{(2s_2-1)\times 1}\\
0&0&A(s_3)&...&0&\underline{1}_{(2s_3-1)\times 1}\\
&&&...
\\0&0&0&...&A(s_k)&\underline{1}_{(2s_k-1)\times 1}
\end{pmatrix}$ is a $2r\times r$ matrix and is irreducible. Then one must have that
$\sum\limits_{i=1}^k(s_i-1)+1=r$ and
$\sum\limits_{i=1}^k(2s_i-1)=2r$. From these two equations, it is
easy to check that $k$ must be 2.

\vskip .2cm The proof of  the last part of Lemma~\ref{ir}(b) is
immediate.
$\hfill\Box$\\

\begin{rem}
By a direct computation, one may  give the first 6 numbers for
$N(r)= |\mathcal{S}_{2r}\backslash\mathcal{M}_r/\text{GL}(r,{\Bbb
Z}_2)|$ as follows:
$$\begin{tabular}{|c|c|c|c|c|c|c|}
 \hline $r$ &
 1 & 2& 3& 4& 5& 6\\
\hline $N(r)$ & 1& 1&  1&  2&  2& 3
 \\
\hline
\end{tabular}$$
  The
representatives of all orbits in
$\mathcal{S}_{2r}\backslash\mathcal{M}_r/\text{GL}(r,{\Bbb Z}_2)$
for $r\leq 6$ are stated as follows: {\tiny
$$\begin{tabular}{|c|c|c|c|c|c|c|}
 \hline $r$ &
 1 & 2& 3& 4& 5&6\\
\hline Rep. & $\begin{pmatrix}1\\1
\end{pmatrix}$& $\begin{pmatrix}I_2\\I_2
\end{pmatrix}$& $\begin{pmatrix}I_3\\I_3
\end{pmatrix}$ & $\begin{pmatrix}I_4\\I_4
\end{pmatrix}$, $B(4)$&$\begin{pmatrix}I_5\\I_5
\end{pmatrix}$,  $\begin{pmatrix}B(4)&0\\
0&\underline{1}_{2\times 1}
\end{pmatrix}$
 &$\begin{pmatrix}I_6\\I_6
\end{pmatrix}$,  $\begin{pmatrix} B(4)&0&0\\
0& \underline{1}_{2\times 1}& 0\\
\\0&0&\underline{1}_{2\times 1}
\end{pmatrix}$, $B(6)$
\\
\hline
\end{tabular}$$}

\noindent where $B(l)=\begin{pmatrix}
A(l)&\underline{1}_{(2l-1)\times 1}\\0&1
\end{pmatrix}$.
Note that for $ r\geq 7$, one also can find the lower bound of
$N(r)$ as follows:
$$\begin{tabular}{|c|c|c|c|c|c|c|c|c|c|}
 \hline $r$ &
 7 &8& 9& 10& 11& 12& 13&14&$\cdots$\\
\hline $N(r)\geq$ & 4&6&  7& 9&  12& 16&20&25&$\cdots$
 \\
\hline
\end{tabular}$$
\end{rem}

\end{document}